\documentclass[a4paper]{article}
\usepackage{amsfonts,amsmath,graphicx,calc}
\usepackage[sort&compress,numbers]{natbib}
\usepackage{amsthm}

\theoremstyle{plain}

\newtheorem{theorem}{Theorem}

\newtheorem{lemma}{Lemma}
\newtheorem{corollary}{Corollary}

\theoremstyle{definition}
\newtheorem{conjecture}{Conjecture}

\newtheorem{open}{Open Problem}

\theoremstyle{remark}

\newcommand{\lemref}[1]{Lemma~\ref{lem:#1}}

\newcommand{\twolemref}[2]{Lemmata~\ref{lem:#1} and \ref{lem:#2}}

\newcommand{\thmref}[1]{Theorem~\ref{thm:#1}}

\newcommand{\secref}[1]{Section~\ref{sec:#1}}

\newcommand{\corref}[1]{Corollary~\ref{cor:#1}}

\newcommand{\eqnref}[1]{\eqref{eqn:#1}}

\newcommand{\seclabel}[1]{\label{sec:#1}}

\newcommand{\thmlabel}[1]{\label{thm:#1}}
\newcommand{\lemlabel}[1]{\label{lem:#1}}

\newcommand{\corlabel}[1]{\label{cor:#1}}
\newcommand{\eqnlabel}[1]{\label{eqn:#1}}

\newcommand{\mySection}[2]{\section{#1}\seclabel{#2}}

\DeclareFontShape{OT1}{cmr}{b}{sc}
{<5><6><7><8><9><10><12><10.95><14.4><17.28><20.74><24.88>cmbcsc10}{}
\DeclareFontShape{OT1}{cmr}{bx}{sc}
{<5><6><7><8><9><10><12><10.95><14.4><17.28><20.74><24.88>cmbcsc10}{}

\DeclareFontShape{OT1}{cmr}{b}{tt}
{<5><6><7><8><9><10><12><10.95><14.4><17.28><20.74><24.88>cmbtt10}{}
\DeclareFontShape{OT1}{cmr}{bx}{tt}
{<5><6><7><8><9><10><12><10.95><14.4><17.28><20.74><24.88>cmbtt10}{}

\newcommand{\NP}{\ensuremath{\mathcal{NP}}}

\newcommand{\Z}{\ensuremath{\mathbb{Z}}}

\newcommand{\FancyFootnotes}{\renewcommand{\thefootnote}{\fnsymbol{footnote}}}

\newcommand{\Prob}[1]{\ensuremath{\textbf{\textup{P}}#1}}

\newcommand{\e}{\ensuremath{\boldsymbol{e}}}

\newcommand{\SET}[1]{\ensuremath{\protect\left\{#1\right\}}}

\newcommand{\bracket}[1]{\ensuremath{\protect\left(#1\right)}}

\newcommand{\ceil}[1]{\ensuremath{\protect\lceil#1\rceil}}
\newcommand{\CEIL}[1]{\ensuremath{\protect\left\lceil#1\right\rceil}}

\newcommand{\half}{\ensuremath{\protect\tfrac{1}{2}}}

\newcommand{\quarter}{\ensuremath{\protect\tfrac{1}{4}}}

\newcommand{\etal}{~et~al.~}

\newcommand{\Oh}[1]{\ensuremath{\protect\mathcal{O}(#1)}}

\newcommand{\Settings}[5]{

\renewcommand{\baselinestretch}{#1}
\setlength{\footnotesep}{\baselinestretch\footnotesep}

\setlength{\marginparwidth}{0mm}
\setlength{\marginparsep}{0mm}
\setlength{\marginparpush}{0mm}

\newlength{\LRmargin}
\setlength{\LRmargin}{#2}
\newlength{\TBmargin}
\setlength{\TBmargin}{#3}

\setlength{\voffset}{-1in}           
\setlength{\hoffset}{-1in}           

\setlength{\headsep}{#4}
\setlength{\footskip}{#5}

\setlength{\oddsidemargin}{\LRmargin}
\setlength{\evensidemargin}{\LRmargin}
\setlength{\topmargin}{\TBmargin}

\setlength{\textheight}{\paperheight}
\addtolength{\textheight}{-1in}			
\addtolength{\textheight}{-1.0\voffset}
\addtolength{\textheight}{-2.0\TBmargin}	
\addtolength{\textheight}{-1.0\headheight}
\addtolength{\textheight}{-1.0\headsep}
\addtolength{\textheight}{-1.0\footskip}

\setlength{\textwidth}{\paperwidth}
\addtolength{\textwidth}{-1.0in}		
\addtolength{\textwidth}{-1.0\hoffset}
\addtolength{\textwidth}{-2.0\LRmargin}		
}

\newlength{\marginboxwidth}

\usepackage{pifont}

\FancyFootnotes
\Settings{1.23}{35mm}{25mm}{0cm}{1cm}

\DeclareMathOperator{\VOL}{vol}
\DeclareMathOperator{\MAG}{mag}

\newcommand{\E}{\mathcal{E}}
\newcommand{\D}{\mathcal{D}}
\renewcommand{\S}{\mathcal{S}}

\begin{document}

\title{\textbf{Drawing a Graph in a Hypercube}}

\author{David~R.~Wood\thanks{
Research supported by the Government of Spain grant MEC SB2003-0270, and by the projects MCYT-FEDER BFM2003-00368 and Gen.\ Cat 2001SGR00224. Research initiated in the Department of Applied Mathematics and the Institute for Theoretical Computer Science at Charles University, Prague, Czech Republic. Supported by project LN00A056 of the Ministry of Education of the Czech Republic, and by the European Union Research Training Network COMBSTRU (Combinatorial Structure of Intractable Problems).}\\
Departament de Matem{\`a}tica Aplicada II\\
Universitat Polit{\`e}cnica de Catalunya\\
Barcelona, Spain\\
\texttt{david.wood@upc.edu}}

\date{November 16, 2004 (revised \today)}
\maketitle

\begin{abstract}    
A \emph{$d$-dimensional hypercube drawing} of a graph represents the vertices by distinct points in  $\{0,1\}^d$, such that the line-segments representing the edges do not cross. We study lower and upper bounds on the minimum number of dimensions in hypercube drawing of a given graph. This parameter turns out to be related to Sidon sets and antimagic injections.
\end{abstract}

\mySection{Introduction}{Introduction}

Two-dimensional graph drawing \citep{DETT99, KaufmannWagner01}, and to a lesser extent, three-dimensional graph drawing \citep{Landgraf01, Wood-PhD, DMW-SJC05} have been widely studied in recent years. Much less is known about graph drawing in higher dimensions. For research in this direction, see references \citep{Wood-PhD, Eppstein-EuJC05, Wood-GD99, BDLS95, EHT65}. This paper studies drawings of graphs in which the vertices are positioned at the points of a hypercube.

We consider undirected, finite, and simple graphs $G$ with vertex set $V(G)$ and edge set $E(G)$. Consider an injection $\lambda:V(G)\rightarrow\{0,1\}^d$.  For each edge $vw\in E(G)$, let $\lambda(vw)$ be the open line-segment with endpoints $\lambda(v)$ and $\lambda(w)$. Two distinct edges $vw,xy\in E(G)$ \emph{cross} if $\lambda(vw)\cap\lambda(xy)\ne\emptyset$. We say $\lambda$ is a \emph{$d$-dimensional hypercube drawing} of $G$ if no two edges of $G$ cross. A $d$-dimensional hypercube drawing is said to have \emph{volume} $2^d$. That is, the volume is the total number of points in the hypercube, and is a measure of the efficiency of the drawing. Let $\VOL(G)$ be the minimum volume of a hypercube drawing of a graph $G$. This paper studies lower and upper bounds on $\VOL(G)$.

The remainder of the paper is organised as follows. In \secref{Antimagic} we review material on Sidon sets and so-called antimagic injections of graphs. In \secref{Basics} we explore the relationship between hypercube drawings and antimagic injections. This enables lower and upper bounds  on $\VOL(K_n)$ to be proved. In \secref{Degeneracy}, we present a simple algorithm for computing an antimagic injection that gives upper bounds on the volume of hypercube drawings in terms of the degeneracy of the graph.  In \secref{Queues} we prove a relationship between antimagic injections and queue layouts of graphs that enables an \NP-completeness result to be concluded. In \secref{Width} we relate antimagic injections of graphs to the bandwidth and pathwidth parameters. Finally, in \secref{Main} we give an asymptotic bound on the volume of hypercube drawings. The proof is based on the Lov{\'a}sz Local Lemma. 

\mySection{Sidon Sets and Antimagic Injections}{Antimagic}


A set $S\subseteq\Z^+$ is called \emph{Sidon} if $a+b=c+d$ implies $\{a,b\}=\{c,d\}$ for all $a,b,c,d\in S$. See the recent survey by \citet{OBryant-EJC04} for results and numerous references on Sidon sets. A graph in which self-loops are allowed (but no parallel edges) is called a \emph{pseudograph}. For a pseudograph $G$, an injection $f:V(G)\rightarrow\mathbb{Z}^+$ is \emph{antimagic} if  $f(v)+f(w)\ne f(x)+f(y)$ for all distinct edges $vw,xy\in E(G)$; see \citep{Gallian03, Wood-AJC02, BolPik-EuJC05}. Let $[k]:=\{1,2,\dots,k\}$. Let $\MAG(G)$ be the minimum $k$ such that the pseudograph $G$ has an antimagic injection $f:V(G)\rightarrow[k]$. 

Let $K_n^+$ be the complete pseudograph; that is, every pair of vertices are adjacent and there is one loop at every vertex. Clearly an antimagic injection of $K_n^+$ is nothing more than a Sidon set of cardinality $n$. It follows from results by \citet{Singer38} and \citet{ET41} (see \citet{BolPik-EuJC05}) that
\begin{equation}
\eqnlabel{CompleteMagic}
\MAG(K_n)\;=\;(1+o(1))n^2\;\;\text{and}\;\;
\MAG(K_n^+)\;=\;(1+o(1))n^2\enspace.
\end{equation}

Note the following simple lower bound.

\begin{lemma}
\lemlabel{MagicLowerBound}
Every pseudograph $G$ satisfies $\MAG(G)\geq\max\{|V(G)|,\half(|E(G)|+3)\}$.
\end{lemma}

\begin{proof} That $\MAG(G)\geq|V(G)|$ follows from the definition.  Let $\lambda:V(G)\rightarrow[k]$ be an antimagic injection of $G$. For every edge $vw\in E(G)$, $\lambda(v)+\lambda(w)$ is a distinct integer in $\{3,4,\dots,2k-1\}$. Thus $|E(G)|\leq2k-3$ and $k\geq\half(|E(G)|+3)$.
\end{proof}

\mySection{Hypercube Drawings}{Basics}

Consider the maximum number of edges in a hypercube drawing. The following observation is a special case of a result by  \citet{BCMW-JGAA04} regarding the volume of grid drawings,  where the bounding box is unrestricted.

\begin{lemma}[\citep{BCMW-JGAA04}]
\lemlabel{Extremal}
The maximum number of edges in a $d$-dimensional hypercube drawing is $3^d-2^d$.
\end{lemma}

Trivially, $\VOL(G)\geq|V(G)|$. For dense graphs, we have the following improved lower bound.

\begin{lemma}
\lemlabel{LowerBound}
Every $n$-vertex $m$-edge graph $G$ satisfies $\VOL(G)\,\geq\,(n+m)^{1/\log_23}\,=\,(n+m)^{0.631\ldots}$.
\end{lemma}

\begin{proof}
Suppose that $G$ has a $d$-dimensional hypercube drawing. By \lemref{Extremal} and since $n\leq 2^d$, we have $n+m\leq 3^d$. That is, $d\geq \log_2(n+m)/\log_23$, and the volume $2^d\geq (n+m)^{1/\log_23}$.
\end{proof}

Now we characterise when two edges cross.

\begin{lemma}
\lemlabel{Crossing}
Consider an injection $\lambda:V(G)\rightarrow\{0,1\}^d$ for some graph $G$. Two distinct edges $vw,xy\in E(G)$ cross if and only if $\lambda(v)+\lambda(w)=\lambda(x)+\lambda(y)$.
\end{lemma}

\begin{proof} Suppose that $\lambda(v)+\lambda(w)=\lambda(x)+\lambda(y)$. Then $\half(\lambda(v)+\lambda(w))=\half(\lambda(x)+\lambda(y))$. That is, the midpoint of $\lambda(vw)$ equals the midpoint of $\lambda(xy)$. Hence $vw$ and $xy$ cross. (Note that this idea is used to prove the upper bound in  \lemref{Extremal}, since the number of midpoints is at most $3^d-2^d$.)\  Conversely, suppose that $vw$ and $xy$ cross. Since all vertex coordinates are $0$ or $1$, the point of intersection between $\lambda(vw)$ and $\lambda(xy)$ is the midpoint of both edges. That is, $\half(\lambda(v)+\lambda(w))=\half(\lambda(x)+\lambda(y))$, and  $\lambda(v)+\lambda(w)=\lambda(x)+\lambda(y)$. 
\end{proof}

Loosely speaking, \lemref{Crossing} implies that a hypercube drawing of $G$ can be thought of as an antimagic injection of $G$ into a set of boolean vectors (where vector addition is \emph{not} modulo $2$).  Moreover, from an antimagic injection we can obtain a hypercube drawing, and vice versa.

\begin{lemma}
\lemlabel{Convert}
Every graph $G$ satisfies $\VOL(G)\leq 2^{\ceil{\log_2\MAG(G)}}<2\MAG(G)$.
\end{lemma}

\begin{proof}
Let $k:=\MAG(G)$, and  let $f:V(G)\rightarrow[k]$ be an antimagic injection of $G$.  For each vertex $v\in V(G)$, let $\lambda(v)$ be the $\ceil{\log_2k}$-bit binary representation of $f(v)$. Suppose that edges $vw$ and $xy$ cross.   By \lemref{Crossing}, $\lambda(v)+\lambda(w)=\lambda(x)+\lambda(y)$.   For each $1\leq i\leq\ceil{\log_2k}$, the sum of the $i$-th coordinates of $v$ and $w$ equals the sum of the $i$-th coordinates of $x$ and $y$.   Thus $f(v)+f(w)=f(x)+f(y)$, which is the desired contradiction.  Therefore no two edges cross, and $\lambda$ is a  \ceil{\log_2k}-dimensional hypercube drawing of $G$.
\end{proof}

\begin{lemma}
\lemlabel{ConvertBack}
Every graph $G$ satisfies $\MAG(G)\leq\VOL(G)^{\log_23}=\VOL(G)^{1.585\ldots}$.
\end{lemma}

\begin{proof} Let $\lambda:V(G)\rightarrow\{0,1\}^d$ be a hypercube drawing of $G$, where $d=\log_2\VOL(G)$. For each vertex $v\in V(G)$, define an integer $f(v)$ so that $\lambda(v)$ is the base-$3$ representation of $f(v)$.  Now $\lambda(v)+\lambda(w)\in\{0,1,2\}^d$. Thus  $\lambda(v)+\lambda(w)=\lambda(x)+\lambda(y)$ if and only if $f(v)+f(w)=f(x)+f(y)$. Since edges do not cross in $\lambda$ and by \lemref{Crossing}, $f$ is an antimagic injection of $G$ into $[3^d]=[3^{\log_2\VOL(G)}]=[\VOL(G)^{\log_23}]$.
\end{proof}

Consider the minimum volume of a hypercube drawing of the complete graph $K_n$.

\begin{lemma}
\lemlabel{BinaryVectors}
Let $V=\{\vec{v}_1,\vec{v}_2,\dots,\vec{v}_n\}$  be a set of binary $d$-dimensional vectors. Then $V$ is the vertex set of a hypercube drawing of $K_n$ if and only if $\vec{v}_i+\vec{v}_j\ne\vec{v}_k+\vec{v}_\ell$ for all distinct pairs $\{i,j\}$ and $\{k,\ell\}$. 
\end{lemma}

\begin{proof}
Suppose that $V$ is the vertex set of a hypercube drawing of $K_n$.  Since no two edges cross, by \lemref{Crossing}, $\vec{v}_i+\vec{v}_j\ne\vec{v}_k+\vec{v}_\ell$ for all distinct pairs $\{i,j\}$ and $\{k,\ell\}$ with $i\ne j$ and $k\ne\ell$. If $i=j$ and $k=\ell$, then $\vec{v}_i+\vec{v}_j\ne\vec{v}_k+\vec{v}_\ell$ because distinct vertices are mapped to distinct points. If $i=j$ and $k\ne\ell$, then $\vec{v}_i+\vec{v}_j\ne\vec{v}_k+\vec{v}_\ell$,  as otherwise the midpoint of the edge $v_kv_\ell$ would coincide with the vertex $v_i$, which is clearly impossible. Hence $\vec{v}_i+\vec{v}_j\ne\vec{v}_k+\vec{v}_\ell$ for all distinct pairs $\{i,j\}$ and $\{k,\ell\}$. The converse result follows immediately from \lemref{Crossing}.
\end{proof}

Sets of binary vectors satisfying \lemref{BinaryVectors}  were first studied by
\citet{Lindstrom-JNT72, Lindstrom-JCT69}, and more recently by
\citet{CLZ-JCTA01}. Their results can be interpreted as follows, where the
lower bound is by \citet{CLZ-JCTA01}, and the upper bound follows from \eqnref{CompleteMagic} and \lemref{Convert}.


\begin{theorem}
\thmlabel{Complete}
Every complete graph $K_n$ satisfies $\VOL(K_n)<(2+o(1))n^2$, and 
$\VOL(K_n)>n^{1.7384\ldots}$ for large enough $n$.\qed 
\end{theorem}

\mySection{Degeneracy}{Degeneracy}

\citet{Wood-AJC02} proved that every  $n$-vertex $m$-edge graph $G$ with maximum degree $\Delta$ satisfies $\MAG(G)<(\Delta(m-\Delta)+n)$. Thus \lemref{Convert} implies that 
\begin{equation}
\eqnlabel{Greedy}
\VOL(G)<2(\Delta(m-\Delta)+n)\enspace.
\end{equation}
This result of \citet{Wood-AJC02} is proved using a greedy algorithm. We can obtain a more precise result as follows. The \emph{degeneracy} of a graph $G$ is the maximum, taken over all induced subgraphs $H$ of $G$, of the minimum degree of $H$.

\begin{lemma}
\lemlabel{Degen}
Every $n$-vertex $m$-edge graph $G$ with degeneracy $d$ satisfies $\MAG(G)\leq n+dm$, and thus $\VOL(G)<2n+2dm$.\qed
\end{lemma}

\begin{proof} We proceed by induction on $n'$ with the hypothesis that ``every induced subgraph $H$ of $G$ on $n'$ vertices has  $\MAG(H)\leq n'+dm$.''  If $n'=1$ the result is trivial.  Let $H$ be an induced subgraph of $G$ on $n'\geq2$ vertices.  Then $H$ has a vertex $v$ of degree at most $d$ in $H$.  By induction,  $H\setminus v$  has an antimagic injection  $\lambda:V(H\setminus v)\rightarrow[n'-1+dm]$. Now
\begin{align*}
&\big|\{\lambda(x):x\in V(H\setminus v)\}\;\cup\;
\{\lambda(x)+\lambda(y)-\lambda(w):xy\in E(H\setminus v),vw\in E(H)\}\big|\\
\leq\;
&|V(H\setminus v)|\,+\,\deg_H(v)\cdot|E(H\setminus v)|\\
\leq\;
&n'-1+dm\enspace.
\end{align*}
Thus there exists an $i\in[n'+dm]$ such that $\lambda(x)\ne i$ for all $x\in V(H\setminus v)$, and $\lambda(x)+\lambda(y)-\lambda(w)\ne i$ for all edges $xy\in E(H\setminus v)$ and $vw\in E(H)$. Let $\lambda(v):=i$. Thus $\lambda(v)\ne\lambda(x)$ for all $x\in V(H)$, and $\lambda(v)+\lambda(w)\ne\lambda(x)+\lambda(y)$ for all edges $xy\in E(H)$ and $vw\in E(G)$. Thus $\lambda$ is an antimagic injection of $H$ into $[n'+dm]$,  and $\MAG(H)\leq n'+dm$. By induction, $\MAG(G)\leq n+dm$.
\end{proof}

Planar graphs $G$ are $5$-degenerate, and thus satisfy $\MAG(G)<16n$ and $\VOL(G)<32n$ by \twolemref{Convert}{Degen}. More generally,  \citet{Kostochka82} and \citet{Thomason84,Thomason01} independently proved that  a graph $G$ with no $K_k$ minor is $\Oh{k\sqrt{\log k}}$-degenerate, and thus satisfy  $\MAG(G)\in\Oh{k^2(\log k)n}$ and $\VOL(G)\in\Oh{k^2(\log k)n}$ by \twolemref{Convert}{Degen}. As we now show, a large clique minor does not necessarily force up $\MAG(G)$ or $\VOL(G)$. Let $K_n'$ be the graph obtained from $K_n$ by subdividing every edge once. Say $K_n'$ has $n':=n+\binom{n}{2}$ vertices. Clearly $K_n'$ is $2$-degenerate. If follows from \lemref{Degen} that $\MAG(K_n')\leq 5n'+o(n')$ and $\VOL(K_n')\leq 10n'+o(n')$, yet $K_n'$ contains a $(\sqrt{2n'}+o(n'))$-clique minor.

\mySection{Queue Layouts and Complexity}{Queues}

Let $G$ be a graph. A bijection  $\sigma:V(G)\rightarrow[|V(G)|]$ is called a \emph{vertex ordering} of $G$. Consider edges $vw,xy\in E(G)$ with no common endpoint.  Without loss of generality $\sigma(v)<\sigma(w)$,  $\sigma(x)<\sigma(y)$ and $\sigma(v)<\sigma(x)$. We say $vw$ and $xy$ are \emph{nested} in $\sigma$ if $\sigma(v)<\sigma(x)<\sigma(y)<\sigma(w)$. A \emph{queue} in $\sigma$ is a set of edges $Q\subseteq E(G)$  such that no two edges in $Q$ are nested in $\sigma$. A \emph{$k$-queue layout} of $G$ consists of a vertex ordering $\sigma$ of $G$, and a partition of $E(G)$ into $k$ queues in $\sigma$. Heath\etal\citep{HLR-SJDM92,HR-SJC92} introduced queue layouts; see \citet{DujWoo-DMTCS04} for references and a summary of known results.

\begin{lemma}
\lemlabel{Queue}
If a graph $G$ has a $1$-queue layout, then $\MAG(G)=|V(G)|$.
\end{lemma}

\begin{proof} Let $\sigma:V(G)\rightarrow[|V(G)|]$ be the vertex ordering in a $1$-queue layout of $G$.  If for distinct edges $vw,xy\in E(G)$, we have $\sigma(v)+\sigma(w)=\sigma(x)+\sigma(y)$, then $vw$ and $xy$ are nested. Since no two edges are nested in a $1$-queue layout, $\sigma$ is an antimagic injection of $G$, and $\MAG(G)\leq|V(G)|$. 
\end{proof}

\citet{HR-SJC92} proved that it is \NP-complete to determine if a given graph has a $1$-queue layout. Thus, \lemref{Queue} implies:

\begin{corollary}
\corlabel{Hard}
Testing whether $\MAG(G)=|V(G)|$ is \NP-complete.\qed
\end{corollary}

It is has been widely conjectured that it is \NP-complete to recognise graphs that admit certain types of magic and antimagic injections. \corref{Hard} is the first result in this direction that we are aware of.

\begin{open}
Every $k$-queue graph $G$ on $n$ vertices is $4k$-degenerate \citep{DujWoo-DMTCS04, Pemmaraju-PhD}. By \lemref{Degen}, $\MAG(G)\in\Oh{k^2n}$ and $\VOL(G)\in\Oh{k^2n}$. Can these bounds be improved to \Oh{kn}?
\end{open}


\mySection{Bandwidth and Pathwidth}{Width}

Let $P_n^k$ be the $k$-th power of a path. That is, $P_n^k$ is the graph with vertex set $\{v_0,v_1,\dots,v_{n-1}\}$ and edge set $\{v_iv_j:1\leq|i-j|\leq k\}$. Now $P_n^k$ has $kn-\half k(k+1)$ edges. By \lemref{MagicLowerBound}, $\MAG(P_n^k)\geq \half(kn-\half k(k+1)+3)$. The following upper bound is a generalisation of the construction of a Sidon set by \citet{ET41}.

\begin{lemma}
\lemlabel{PathPower}
For every prime $p$, $\MAG(P_n^p)\leq p(2n-1)$.
\end{lemma}

\begin{proof}
If $p=2$ then $\MAG(P_n^2)$ has a $1$-queue layout, and $\MAG(P_n^2)=n$ by \lemref{Queue}. Now assume that $p>2$.  Let $\lambda(v_i):=1+2pi+(i^2\bmod{p})$ for every vertex $v_i$, $0\leq i\leq n-1$. Clearly $\lambda$ is an injection into $[p(2n-1)]$. Suppose on the contrary, that there are distinct edges $v_iv_\ell$ and $v_jv_k$ with  $\lambda(v_i)+\lambda(v_\ell)=\lambda(v_j)+\lambda(v_k)$. Without loss of generality, $i<j<k<\ell\leq i+p$. Then
\begin{equation*}
2pi+(i^2\bmod{p})+2p\ell+(\ell^2\bmod{p})=2pj+(j^2\bmod{p})+2pk+(k^2\bmod{p})\enspace.
\end{equation*}
That is,
\begin{equation*}
2p(i+\ell-j-k)=(j^2\bmod{p})+(k^2\bmod{p})-(i^2\bmod{p})-(\ell^2\bmod{p})\enspace.
\end{equation*}
Now  $|(j^2\bmod{p})+(k^2\bmod{p})-(i^2\bmod{p})-(\ell^2\bmod{p})|\leq2(p-1)$.
Thus $i+\ell-j-k=0$, and
\begin{equation*}
(i^2\bmod{p})+(\ell^2\bmod{p})=(j^2\bmod{p})+(k^2\bmod{p})\enspace.
\end{equation*}
Thus
\begin{equation}
\eqnlabel{PPP}
i^2+\ell^2\equiv j^2+k^2\pmod{p}\enspace.
\end{equation}
Let $a:=j-i$ and $b:=k-i$. Then $0<a<b<p$. Since $i+\ell=j+k$, we have $\ell=i+a+b$. Rewriting \eqnref{PPP},
\begin{equation*}
i^2+(i+a+b)^2\equiv(i+a)^2+(i+b)^2\pmod{p}\enspace.
\end{equation*}
Hence $2ab\equiv0\pmod{p}$. Since $p$ is prime and $p>2$,  $a\equiv0\pmod{p}$ or $b\equiv0\pmod{p}$,  which is a contradiction since $0<a<b<p$. Hence $\lambda(v_i)+\lambda(v_\ell)\ne\lambda(v_j)+\lambda(v_k)$, and $\lambda$ is antimagic.
\end{proof}

The \emph{bandwidth} of an $n$-vertex graph $G$ is the minimum $k$ such that $G$ is a subgraph of $P_n^k$. By Bertrand's postulate there is a prime $p\leq2k$. Thus \twolemref{Convert}{PathPower} imply:

\begin{corollary}
\corlabel{Bandwidth}
Every $n$-vertex graph $G$ with bandwidth $k$ has $\MAG(G)\leq2k(2n-1)$ and $\VOL(G)<4k(2n-1)$.\qed
\end{corollary}

We have the following technical lemma.

\begin{lemma}
\lemlabel{Technical}
Let $G$ be a graph.  Let $f_V:V(G)\rightarrow[t]\times[r]$ be an injection. Define a function $f_E:E(G)\rightarrow\binom{[t]}{2}\times[2r]$ as follows. For every edge $vw\in E(G)$ with $f_V(v)=(a,i)$ and $f_V(w)=(b,j)$, let $f_E(vw):=(\{a,b\},i+j)$. If $f_E$ is also an injection, then $\MAG(G)\leq (2+o(1))t^2r$.
\end{lemma}

\begin{proof}
\citet{Singer38} proved that there is a Sidon set $\{s_1,s_2,\dots,s_t\}\in[(1+o(1))t^2]$. For every vertex $v\in V(G)$ with $f(v)=(a,i)$,  let $\lambda(v):=2r(s_a-1)+i$. Since $f$ is an injection, $\lambda$ is an injection into $[(2+o(1))t^2r]$. We claim that $\lambda$ is antimagic. Suppose on the contrary that there are distinct edges $vw,xy\in E(G)$ with $\lambda(v)+\lambda(w)=\lambda(x)+\lambda(y)$. Say $f(v)=(a,i)$, $f(w)=(b,j)$, $f(x)=(c,k)$, and $f(y)=(d,\ell)$. Then
\begin{equation}
\eqnlabel{SomeEqn}
2r(s_a-1)+i\,+\,
2r(s_b-1)+j\,=\,
2r(s_c-1)+k\,+\,
2r(s_d-1)+\ell\enspace.
\end{equation}
That is, $2r(s_a+s_b-s_c-s_d)\,=\,k+\ell-i-j$. Now $|k+\ell-i-j|<2r$. Thus $s_a+s_b=s_c+s_d$. Since  $\{s_1,s_2,\dots,s_t\}$ is Sidon, $\{a,b\}=\{c,d\}$. By \eqnref{SomeEqn}, $i+j=k+\ell$. Hence, $f_E(vw)=f_E(xy)$, which is a contradiction since $f_E$ is an injection by assumption. Thus $\lambda(v)+\lambda(w)\ne\lambda(x)+\lambda(y)$, and $\lambda$ is antimagic.  Hence $\MAG(G)\leq(2+o(1))t^2r$.
\end{proof}

Let $\S$ be a set of closed intervals in $\mathbb{R}$. Associated with $\S$,
is the \emph{interval graph} with vertex set $\S$ such that two vertices are
adjacent if and only if the corresponding intervals have a non-empty intersection. The \emph{pathwidth} of a graph $G$ is the minimum $k$ such that $G$ is a spanning subgraph of an interval graph with no clique on $k+2$ vertices.

\begin{theorem}
Every $n$-vertex graph $G$ with pathwidth $k$ satisfies $\MAG(G)\leq(8+o(1))kn$ and $\VOL(G)\leq(16+o(1))kn$. For all $k$ and $n\geq k+1$, there exist $n$-vertex graphs $G$ with pathwidth $k$ and $\MAG(G)\geq\half kn -\Oh{k^2}$.
\end{theorem}

\begin{proof}
\citet{DMW-SJC05} proved that there is an injection $f$ satisfying  \lemref{Technical} with $t=2k+2$ and $r=\ceil{n/k}$. In fact, they proved the stronger result that for all edges $vw,xy\in E(G)$ with $f(v)=(a,i)$, $f(w)=(b,j)$, $f(x)=(a,k)$, $f(y)=(b,\ell)$, if $i<k$ then $j\leq \ell$ (which implies that $i+j<k+\ell$). By \lemref{Technical}, $\MAG(G)\leq(2+o(1))(2k+2)^2r=(8+o(1))kn$. By \lemref{Convert}, $\VOL(G)\leq(16+o(1))kn$. For the lower bound, let $G=P_n^k$ for example. Then $G$ has pathwidth $k$ and $kn-\half k(k+1)$ edges. By \lemref{MagicLowerBound}, $\MAG(G)\geq\half kn -\Oh{k^2}$.
\end{proof}

\begin{open}
\lemref{Degen} implies  that graphs $G$ of treewidth $k$
satisfy $\MAG(G)\in\Oh{k^2n}$ and $\VOL(G)\in\Oh{k^2n}$. Can these bounds be improved to \Oh{kn}? Note that \citet{Wood-AJC02} proved that every tree $G$ satisfies $\MAG(G)=|V(G)|$, which implies that $\VOL(G)<2|V(G)|$ by \lemref{Convert}. 
\end{open}

\mySection{An Asymptotic Upper Bound}{Main}

Our upper bounds on $\VOL(G)$ have thus far been obtained as corollaries of
upper bounds on $\MAG(G)$. The next theorem, which improves upon 
\eqnref{Greedy}, only applies to hypercube drawings. In fact,  the method used
only gives a \Oh{n+\Delta m} bound on $\MAG(G)$.

\begin{theorem}
\thmlabel{Main}
Every $n$-vertex $m$-edge graph $G$ with maximum degree $\Delta$ satisfies
\begin{equation*}
\VOL(G)\;\leq\;
\Oh{n+(\Delta m)^{1/\log_28/3}}\;=\;
\Oh{n+(\Delta m)^{0.707\ldots}}
\enspace.
\end{equation*}
\end{theorem}

\thmref{Main} is proved using the Local Lemma of \citet{EL75} (see
\cite{MR02}). 

\begin{lemma}[\citep{EL75}]
\lemlabel{LLL}
Let $\E=\{A_1,A_2,\dots,A_n\}$ be a set of `bad' events in some probability
space, such that each event $A_i$ is mutually independent of
$\E\setminus(\{A_i\}\cup\D_i)$ for some $\D_i\subseteq \E$. 
Suppose there is a set $\{x_i\in[0,1):1\leq i\leq n\}$,  such that for all $i$,
\begin{equation}
\eqnlabel{LLL}
\Prob(A_i)\;\leq\;x_i\cdot\!\!\!\!\prod_{A_j\in\D_i}\!\!\!\!(1-x_j)\enspace.
\end{equation}
Then 
\begin{equation*}
\Prob\bracket{\bigwedge_{i=1}^n \overline{A_i}}\;\geq\;\prod_{i=1}^n(1-x_i)\;>\;0\enspace.
\end{equation*}
That is, with positive probability, no event in $\mathcal{E}$ occurs.
\end{lemma}

\begin{proof}[Proof of \thmref{Main}]  Let $d$ be a
positive integer, to be specified later. For each vertex $v\in V(G)$, let
$\lambda(v)$ be a point in $\{0,1\}^d$ chosen randomly and independently.
(One can think of this process as $d$ fair coin tosses for each vertex.)\  
We now set up an application of  \lemref{LLL}. For all pairs of distinct
vertices $v,w\in V(G)$, let $A_{v,w}$ be the event that
$\lambda(v)=\lambda(w)$.  For all disjoint edges $vw,xy\in E(G)$, let
$B_{vw,xy}$ be the event that  $vw$ and $xy$ cross.  

We will apply \lemref{LLL} to prove that with positive probability, no
event occurs. Hence there exists $\lambda$ such that no event occurs. No
$A$-event means that $\lambda$ is an injection. No $B$-event means that
no edges cross. Thus $\lambda$ is a $d$-dimensional hypercube drawing. 

Observe that $\Prob(A_{v,w})=(\half)^d$.
It is easily seen that $\Prob(B_{vw,xy})\leq(\half)^d$.
Below we prove that $\Prob(B_{vw,xy})=(\frac{3}{8})^d$.
The idea here is that it is unlikely that some edges are involved in a crossing.
For example, the actual edges of the hypercube cannot be in a crossing.

Let $M:=\{(x_1,x_2,\dots,x_d):x_i\in\{0,1,2\},1\leq i\leq d\}$. 
Consider an edge $vw\in E(G)$. Clearly $\lambda(v)+\lambda(w)\in M$.
The $i$-coordinate of $\lambda(v)+\lambda(w)$
equals $1$ if and only if the $i$-coordinates of $\lambda(v)$ and $\lambda(w)$
are distinct, which occurs with probability $\half$. 
The $i$-coordinate of $\lambda(v)+\lambda(w)$ equals $0$ if and only if
the $i$-coordinates of $\lambda(v)$ and $\lambda(w)$ both equal $0$, 
which occurs with probability $\quarter$. 
The $i$-coordinate of $\lambda(v)+\lambda(w)$ equals $2$ if and only if
the $i$-coordinates of $\lambda(v)$ and $\lambda(w)$ both equal $1$, 
which occurs with probability $\quarter$. 

Let $M_k$ be the subset of $M$ consisting of those points  with exactly $k$
coordinates equal to $1$. Thus, for every edge $vw\in E(G)$ and point $p\in
M_k$,
\begin{equation*}
\Prob(\lambda(v)+\lambda(w)=p)\;=\;(\half)^k(\quarter)^{d-k}\;=\;2^{k-2d}\enspace.
\end{equation*}
Hence for all disjoint edges $vw,xy\in E(G)$ and points $p\in M_k$,
\begin{equation*}
\Prob(\lambda(v)+\lambda(w)=\lambda(x)+\lambda(y)=p)\;=\;
2^{2k-4d}\enspace.
\end{equation*}
Now $|M_k|=\binom{d}{k}2^{d-k}$. Thus,
\begin{equation*}
\Prob(\lambda(v)+\lambda(w)=\lambda(x)+\lambda(y)\in M_k)\;=\;
\binom{d}{k}2^{d-k}\cdot 2^{2k-4d}\;=\;
\binom{d}{k}2^{k-3d}
\enspace.
\end{equation*}
Thus by \lemref{Crossing},
\begin{equation*}
\Prob(B_{vw,xy})
\;=\;\Prob(\lambda(v)+\lambda(w)=\lambda(x)+\lambda(y))
\;=\;\sum_{k=0}^d\binom{d}{k}2^{k-3d}
\;=\;\bracket{\frac{3}{8}}^d
\enspace.
\end{equation*}

The base of the natural logarithm $\e$ satisfies
the following well-known inequality for all $y>0$:
\begin{equation}
\eqnlabel{eee}
\tfrac{1}{\e}\;<\;\bracket{1-\tfrac{1}{y+1}}^y\enspace.
\end{equation}
Now define
\begin{equation}
\eqnlabel{d}
d\;:=\;\CEIL{\max\SET{
\log_2\e(4n+1),\,
\log_{8/3}\e^2(4\Delta m+1)}}\enspace.
\end{equation}

For each $A$-event, let $x_A:=1/(4n+1)$. 
For each $B$-event, let $x_B:=1/(4\Delta m+1)$.
Thus $0<x_A<1$ and $0<x_B<1$, as required.

Each vertex is involved in at most $n$ $A$-events, and at most $\Delta m$
$B$-events. An $A$-event involves two vertices, and is thus  dependent on at
most $2n$ other $A$-events, and at most $2\Delta m$ $B$-events. A
$B$-event involves four vertices, and is thus dependent on at most $4n$
$A$-events, and on at most $4\Delta m$ other $B$-events. 
We first verify \eqnref{LLL} for each event $A_{v,w}$. By \eqnref{eee}, \begin{equation*}
x_A\bracket{1-x_A}^{2n}\bracket{1-x_B}^{2\Delta m}
\;=\;\frac{1}{4n+1}\bracket{1-\frac{1}{4n+1}}^{2n}\bracket{1-\frac{1}{4\Delta m+1}}^{2\Delta m}
\;\geq\;\frac{1}{\e(4n+1)}\enspace.
\end{equation*}
By the definition of $d$ in \eqnref{d}, $\frac{1}{\e(4n+1)}\geq\frac{1}{2^d}$, and thus
\begin{equation*}
x_A\bracket{1-x_A}^{2n}\bracket{1-x_B}^{2\Delta m}
\;\geq\;\bracket{\frac{1}{2}}^d
\;=\;\Prob(A_{v,w})\enspace.
\end{equation*}
Now we verify \eqnref{LLL} for each event $B_{vw,xy}$. By \eqnref{eee},
\begin{equation*}
      x_B\bracket{1-x_A}^{4n}\bracket{1-x_B}^{4\Delta m}
\;=\;\frac{1}{4\Delta m+1}\bracket{1-\frac{1}{4n+1}}^{4n}\bracket{1-\frac{1}{4\Delta m+1}}^{4\Delta m}
\;\geq\;\frac{1}{\e^2(4\Delta m+1)}\enspace.
\end{equation*}
Note that \eqnref{d} implies that
$\bracket{\tfrac{8}{3}}^d\geq
\e^2(4\Delta m+1)$. Thus,
\begin{equation*}
x_B\bracket{1-x_A}^{4n}\bracket{1-x_B}^{4\Delta m}
\;\geq\;\bracket{\frac{3}{8}}^d
\;=\;\Prob(B_{vw,xy})\enspace.
\end{equation*}
By \lemref{LLL}, there is a $d$-dimensional hypercube drawing of $G$.
The volume $2^d$
is $\Oh{n+(\Delta m)^{1/\log_2{8/3}}}$.
This completes the proof of \thmref{Main}.
\end{proof}






\def\soft#1{\leavevmode\setbox0=\hbox{h}\dimen7=\ht0\advance \dimen7
  by-1ex\relax\if t#1\relax\rlap{\raise.6\dimen7
  \hbox{\kern.3ex\char'47}}#1\relax\else\if T#1\relax
  \rlap{\raise.5\dimen7\hbox{\kern1.3ex\char'47}}#1\relax \else\if
  d#1\relax\rlap{\raise.5\dimen7\hbox{\kern.9ex \char'47}}#1\relax\else\if
  D#1\relax\rlap{\raise.5\dimen7 \hbox{\kern1.4ex\char'47}}#1\relax\else\if
  l#1\relax \rlap{\raise.5\dimen7\hbox{\kern.4ex\char'47}}#1\relax \else\if
  L#1\relax\rlap{\raise.5\dimen7\hbox{\kern.7ex
  \char'47}}#1\relax\else\message{accent \string\soft \space #1 not
  defined!}#1\relax\fi\fi\fi\fi\fi\fi} \def\cprime{$'$}

\end{document}